\documentclass{article}

\usepackage{graphicx}
\usepackage{amsmath}
\usepackage{amsfonts}
\usepackage{amssymb}
\usepackage{float}

\usepackage{color}

\usepackage{cleveref}
\usepackage{epsfig} %% for loading postscript figures
\usepackage{float}
\usepackage{subcaption}
\title{On the Stability of Fractal Differential Equations}
\author{Cemil Tunc $^{1}$,  Alireza Khalili Golmankhaneh $^{2}$  \footnote{Corresponding Author: alirezakhalili2002@yahoo.co.in } \\
  $^{1}$ Department of Mathematics,\\ Faculty of Sciences, Van Yuzuncu Yil University, 65080, Turkey\\
  $^{2}$ Department of Physics, \\ Urmia Branch, Islamic Azad University, Urmia, Iran
}

\begin{document}
\maketitle

\abstract{In this paper, we give a review of fractal calculus which
is an expansion of standard calculus. Fractal calculus is applied for functions which are not differentiable or integrable on totally
disconnected fractal sets such as middle-$\mu$ Cantor sets. Analogues of the Lyapunov functions and features are given for asymptotic behaviors of fractal differential equations. Stability of fractal differentials in the sense of Lyapunov is defined. For the suggested fractal differential equations, sufficient conditions for the stability and uniform boundedness and convergence of the solutions are presented and proved. We present examples and graphs for more details of the results.}

\textbf{Keywords:} Fractal Calculus, Staircase function,Cantor-like sets,Fractal Stability, Fractal Convergence\\
\textbf{MSC[2010]:} 28A80,81Q35, 28A78, 35B40,35B35

\section{Introduction}
Fractals are fragmented shapes at all scales with self-similarity properties theirs fractal dimension exceeds their topological dimension \cite{Mandelbrot,Barnsley,Cattani}. Fractals appear in chaotic dynamical systems as the attractors \cite{Mohammad}. The global attractors of porous media equations and its fractal dimension which is finite under some conditions were suggested in \cite{Efendiev-1,Efendiev,Efendiev-3,Efendiev-5}. The distance of pre-fractal and fractal was derived in terms of the preselected parameters \cite{Efendiev-4}.

Non-standard analysis used to build  the curvilinear coordinate
along the fractal curves (i.e. Ces\`{a}ro and Koch curves) \cite{Nottale}. Theory of scale relativity  suggests  the quantum mechanics formalism as mechanics for fractal space-time \cite{Nottale-1}. Analysis on fractals was studied by using probability theory, measure theory, harmonic analysis, and fractional spaces
\cite{Falconer1,Falconer2,Edgar2,Kigami1,Stillinger,Tarasov}.

Using fractional calculus, electromagnetic fields were provided for fractal charges as generalized distributions and applied to different branches of physics with fractal structures \cite{Zubair,Chen-1,Chen-2}. Non-local fractional derivatives do not have any geometrical and physical meaning so far \cite{Herrmann,Das}.

Local fractional derivatives are needed in many physical problems. The effort of defining local fractional calculus leads to new a measure on fractals \cite{Kolwankar-1,Kolwankar-2}.

In view of this new measure, the $C^{\mu}$-Calculus was formulated for totally disconnected fractal sets and non-differentiable fractal curves \cite{Gangal-1,Gangal-2,Gangal-3,Alireza-Fernandez-1}
During the last decade, several researchers have explored in this area and applied it in different branches of science and engineering
\cite{Golmankhaneh-Paper-1,Golmankhaneh-Paper-2}.~Fractal differential equations (FDEs) were solved and analogous existence and uniqueness theorems were suggested and proved \cite{Golmankhaneh-2s,Golmankhaneh-3s}. The stability of the  impulsive and Lyapunov functions in the sense of  Riemann-Liouville like fractional difference equations were studied in \cite{Wu-1,Wu-3,Hristova,Khaliq}

Motivated by the works mentioned above, we give analogues of asymptotic behaviors of the solutions of FDEs. The stability and asymptotic behavior of differential equations have an important role in various applications in science and engineering. The Lyapunov's second method was applied to show uniform boundedness and convergence to zero of all solutions of second-order non-linear  differential equation \cite{Cemil-1,Cemil-2}. The reader is advised to see the references cited in \cite{Yoshizawa,Hara}.\\
Our aim in this work is to give sufficient conditions for the solutions of FDEs to be uniformly bounded and for the solutions with fractal derivatives to go to zero as $t\rightarrow\infty$. \\
The outline of the manuscript is as follows:\\
In Section 2 we give basic tools and definition we need in the paper.
In Section 3 we define fractal Lyapunov stability and function.
Section 4 gives asymptotic behaviors and conditions for the solutions of FDEs. We present the conclusion of the paper in Section 5.
\section{Preliminaries}
In this section, fractal calculus is summarized which is called   generalized Riemann calculus \cite{Gangal-1,Gangal-2,Gangal-3,Alireza-Fernandez-1}. Fractal calculus  expands standard calculus to involve functions with totally disconnected fractal sets and non-differentiable curves such as Koch and Ces\`{a}ro curves. Fractal calculus was applied for the function with Cantor-like sets with zero Lebesgue measures and non-integer Hausdorff dimensions \cite{Alireza-Fernandez-1,Robert}.
\subsection{The Middle-$\mu$ Cantor set}
The Cantor-like sets contain totally disconnected sets such as  thin fractal, fat fractal, Smith-Volterra-Cantor, k-adic-type, and rescaling Cantor sets \cite{Robert}. The  middle-$\mu$ Cantor set is  obtained by following process \cite{Robert}:\\
First,  delete  an open interval of length $0<\mu<1$ from the middle of the $I=[0,1]$.\\
\begin{equation*}\label{nb}
  C^{\mu}_{1}=[0,\frac{1}{2}(1-\mu)]\cup[\frac{1}{2}(1+\mu),1].
\end{equation*}
Second, remove two  disjoint  open intervals of length $\mu$ from the middle of the remaining closed intervals of first step.\\
\begin{eqnarray*}\label{mnj}
   C^{\mu}_{2}&=&[0,\frac{1}{4}(1-\mu)^2]\cup[\frac{1}{4}(1-\mu)(1+\mu),
   \frac{1}{2}(1-\mu)]\\
   &&\cup[\frac{1}{2}(1+\mu),\frac{1}{2}
   ((1+\mu)+\frac{1}{2}(1+\mu)^2)]\\
   &&\cup[\frac{1}{2}(1+\mu)(1+\frac{1}{2}(1-\mu)),1]
\end{eqnarray*}
\hspace{3cm}~~~~~~~~ \vdots.\\
In $m^{th}$ stage, omit $2^{m-1}$ disjoint open intervals of length $\mu$ from the middle of the remaining closed intervals (See Figure \ref{1}).\\
Finally, we have middle-$\mu$ Cantor set as follows:\\
\begin{equation*}
  C^{\mu}=\bigcap_{m=1}^{\infty}C^{\mu}_{m}.
\end{equation*}
The \textbf{Lebesgue measure} of  $C^{\mu}$ set is  given by
\begin{eqnarray*}\label{nyhu985}
  m(C^{\mu})&=&1-\mu-2(\frac{1}{2}(1-\mu)\mu)-4(\frac{1}{4}(1-\mu)^2\mu)-...\\
  &=&1-\mu\frac{1}{1-(1-\mu)}=1-1=0.
\end{eqnarray*}
The \textbf{Hausdorff dimension} of middle-$\mu$ Cantor set using Hausdorff measure is given by
\begin{equation*}\label{b9q}
 D_{H}(C^{\mu})=\frac{\log(2)}{\log(2)-\log(1-\beta)},
\end{equation*}
where  $H$  indicates Hausdorff measure \cite{Falconer1,Falconer2,Robert}.\\
\subsection{Local fractal calculus}
\begin{figure}[H]
\centering	
	\begin{subfigure}[t]{0.4\textwidth}
		\centering
		\includegraphics[width=\textwidth]{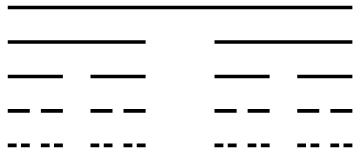}
		\caption{ Middle-$\mu$ Cantor set with $\mu=1/5$}\label{1}
\end{subfigure}
	\quad
\begin{subfigure}[t]{0.4\textwidth}
	\centering
	\includegraphics[width=\textwidth]{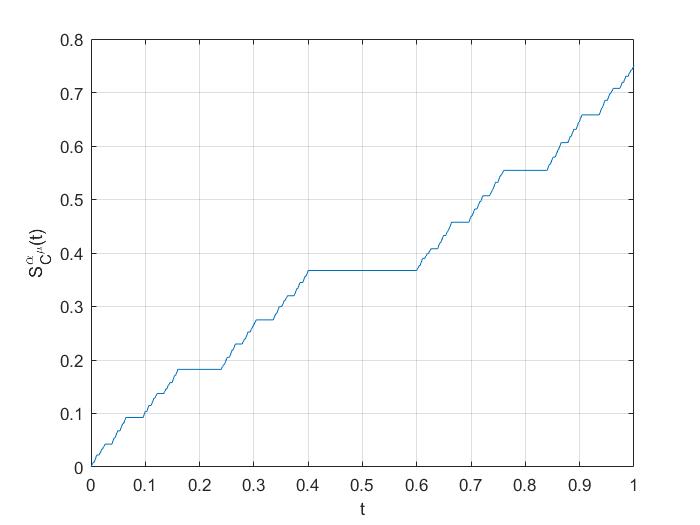}
	\caption{Staircase function corresponding to middle-$\mu$ Cantor set with $\mu=1/5$}\label{2}
	\end{subfigure}
\quad
	\begin{subfigure}[t]{0.4\textwidth}
		\centering
		\includegraphics[width=\textwidth]{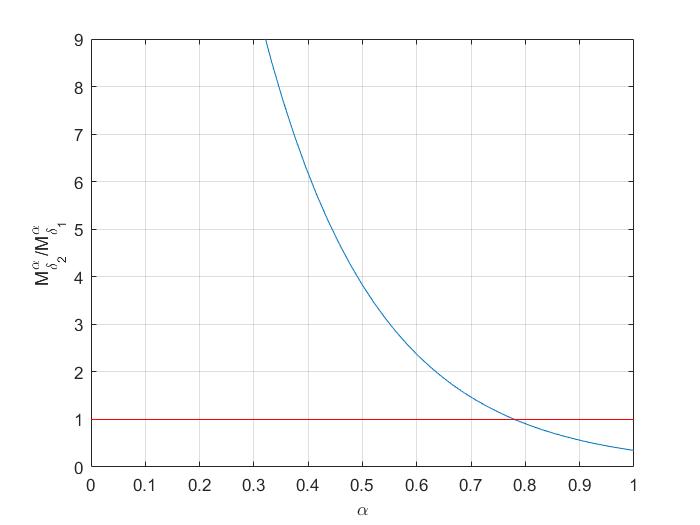}
		\caption{ The $\gamma$-dimension gives $\alpha= 0.75$ to middle-$\mu$ Cantor set with $\mu=1/5$}\label{474}
	\end{subfigure}
\quad
	\begin{subfigure}[t]{0.4\textwidth}
	\centering
		\includegraphics[width=\textwidth]{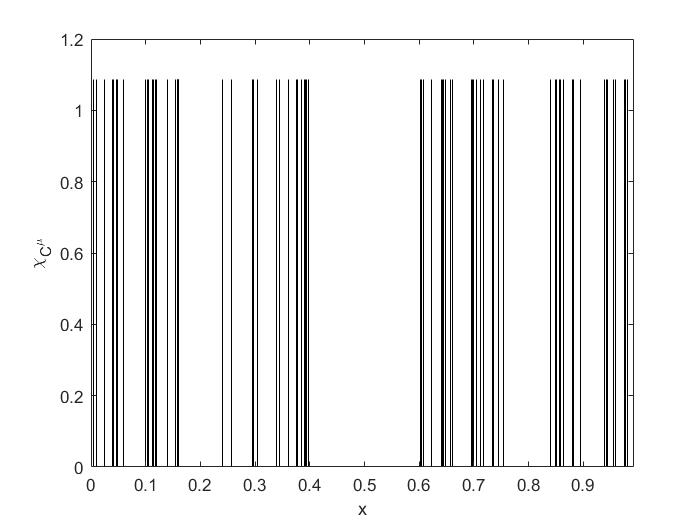}
		\caption{Characteristic function for middle-$\mu$ Cantor set with $\mu=1/5$ }\label{4}
	\end{subfigure}
\caption{ Graphs corresponding to middle-$\mu$ Cantor set with $\mu=1/5$ }\label{poi}
\end{figure}

 The \textbf{flag function} of $C^{\mu}$ is defined by \cite{Gangal-1,Gangal-2},
\begin{equation*}
F(C^{\mu},J)=\begin{cases}
    1 ~~~\textmd{if}~~~ C^{\mu}\cap J \neq \emptyset\\
    0 ~~~~\textmd{otherwise},
\end{cases}
\end{equation*}
where $J=[c_{1},c_{2}]$.~Let $Q_{[c_{1},c_{2}]}=\{c_{1}=t_0,t_1,t_2,\dots,t_m=c_{2}\}$ be a subdivision of $J$. Then, $\mathbb{L}^{\alpha}[C^{\mu},Q]$ is  defined in \cite{Gangal-1,Gangal-2,Alireza-Fernandez-1} by
\begin{equation*}
\mathbb{L}^{\alpha}[C^{\mu},Q]=\sum_{i=1}^{m}\Gamma
(\alpha+1)
(t_{i}-t_{i-1})^{\alpha}F(C^{\mu},[t_{i-1},t_{i}]),
\end{equation*}
where $0<\alpha\leq 1$. \\
The \textbf{mass function} of $C^{\mu}$ is defined in \cite{Gangal-1,Gangal-2,Alireza-Fernandez-1} by
\begin{eqnarray*}\label{s}
    \mathcal{M}^{\alpha}(C^{\mu},c_{1},c_{2})&=&\lim_{\delta\rightarrow0}
    \mathcal{M}^{\alpha}_{\delta}(C^{\mu},c_{1},c_{2})\nonumber\\
    &=&\lim_{\delta\rightarrow0} \left(\inf_{Q_{[c_{1},c_{2}]}:
|Q|\leq\delta}\mathbb{L}^{\alpha}[C^{\mu},Q]\right),
\end{eqnarray*}
here, infimum is taking over all subdivisions $Q$ of $[c_{1},c_{2}]$ satisfying $|Q|:=\max_{1\leq i\leq m}(t_{i}-t_{i-1})\leq\delta$.\\
The \textbf{integral staircase function}  of $C^{\mu}$ is defined in \cite{Gangal-1,Gangal-2,Alireza-Fernandez-1} by
\begin{equation*}\label{t}
    S_{C^{\mu}}^{\alpha}(t)=\begin{cases}
    \mathcal{M}^{\alpha}(C^{\mu},t_{0},t) ~~~~\text{if} ~~~~~t\geq t_{0}\\
    -\mathcal{M}^{\alpha}(C^{\mu},t_{0},t) ~~~~\text{otherwise},
\end{cases}
\end{equation*}
where $t_{0}$ is an arbitrary and fixed real number (See Figure \ref{2}).\\
The \textbf{$\gamma$-dimension} of $C^{\mu}\cap[c_{1},c_{2}]$ is
\begin{align*}%\label{45}
    \dim_{\gamma}(C^{\mu}\cap[c_{1},c_{2}])&=\inf\{\alpha:
    \mathcal{M}^{\alpha}(C^{\mu},c_{1},c_{2})=0\}\nonumber\\
    &=\sup\{\alpha:\mathcal{M}^{\alpha}(C^{\mu},c_{1},c_{2})=\infty\}.
\end{align*}
Figure \ref{474} presents  approximate $\mathcal{M}^{\alpha}_{\delta_{2}}/\mathcal{M}^{\alpha}_{\delta_{1}}$, where $\delta_{2}<\delta_{1}$.~This gives us   $\gamma$-dimension since that value converging to the finite number as $\delta\rightarrow 0$. This result  can also be concluded by choosing different various pairs of $(\delta_{1},\delta_{2})$.\\
The \textbf{ characteristic function}  $\chi_{C^{\mu}}(\alpha,t):\Re\rightarrow \Re  $ is defined  by
\begin{equation*}\label{nhy2jj}
  \chi_{C^{\mu}}(\alpha,t)=\left\{
                             \begin{array}{ll}
                               \frac{1}{\Gamma(\alpha+1)}, & t\in C^{\mu}; \\
                               0, & otherwise.
                             \end{array}
                           \right.
  .
\end{equation*}
In Figure \ref{4} we have ploted the  characteristic function for the middle-$\mu$ choosing $\mu=1/5$.\\
The \textbf{$C^{\alpha}$-limit } of a function $h:\Re \rightarrow \Re$ as $z\rightarrow t$ is defined in \cite{Gangal-1,Gangal-2,Alireza-Fernandez-1} by
\begin{equation*}\label{aw}
    z\in C^{\mu},~~~\forall~ \epsilon,~~~\exists~ \delta,~~~ |z-t|<\delta\Rightarrow|h(z)-l|<\epsilon.
\end{equation*}
If $l$ exists, then we can write
\begin{equation*}\label{we}
    l=\underset{z\rightarrow t}{\operatorname{C^{\mu}_{-}\lim}}~h(z).
\end{equation*}
The \textbf{$C^{\mu}$-continuity} of a  function $h$  at $t\in C^{\mu}$  is  defined in \cite{Gangal-1,Gangal-2}  by
\begin{equation*}
    h(t)=\underset{z\rightarrow t}{\operatorname{C^{\mu}_{-}\lim}}~ h(z).
\end{equation*}
The \textbf{$C^{\mu}$-Differentiation} of a function $h$ at $t\in C^{\mu}$ is defined in \cite{Gangal-1,Gangal-2,Alireza-Fernandez-1} by
\begin{equation*}\label{nb258}
    D_{C^{\mu}}^{\alpha}h(t)=\begin{cases}
    \underset{z\rightarrow t}{\operatorname{C^{\mu}_{-}\lim}} \frac{h(z)-h(t)}{S_{C^{\mu}}^{\alpha}(z)-S_{C^{\mu}}^
    {\alpha}(t)}, ~~~\text{if }z\in C^{\mu},\\
    0, ~~~~~~~~~~~~~~~~~~~~~~~~~~~~~~~\textmd{otherwise},
    \end{cases}
 \end{equation*}
 if  limit exists.\\
The \textbf{ $C^{\mu}$-integral} of $h$ on $[c_{1},c_{2}]$ is  denoted by $\int_{c_{1}}^{c_{2}}h(t) d_{C^{\mu}}^{\alpha}t$ and is approximately  given in \cite{Gangal-1,Gangal-2,Alireza-Fernandez-1} by
\begin{align*}\label{cft1}
 \int_{c_{1}}^{c_{2}}h(t) d_{C^{\mu}}^{\alpha}t \approx \sum_{i=1}^{m} h_{i}(t)(S_{C^{\mu}}^{\alpha}(t_{j})-S_{C^{\mu}}^{\alpha}(t_{j-1})).
\end{align*}
We refer the reader for more meticulous definitions to see in \cite{Gangal-1,Gangal-2,Alireza-Fernandez-1}.
In Figure \ref{poi} we have sketched the middle-$\mu$ Cantor, the staircase function, the characteristic function, and graph of $\mathcal{M}^{\alpha}_{\delta_{2}}/\mathcal{M}^{\alpha}_{\delta_{1}}$ versus to $\alpha$.\\
\section{Fractal Lyapunov stability} 
\begin{figure}[H]
\centering	
\begin{subfigure}[t]{0.4\textwidth}
		\centering
		\includegraphics[width=\textwidth]{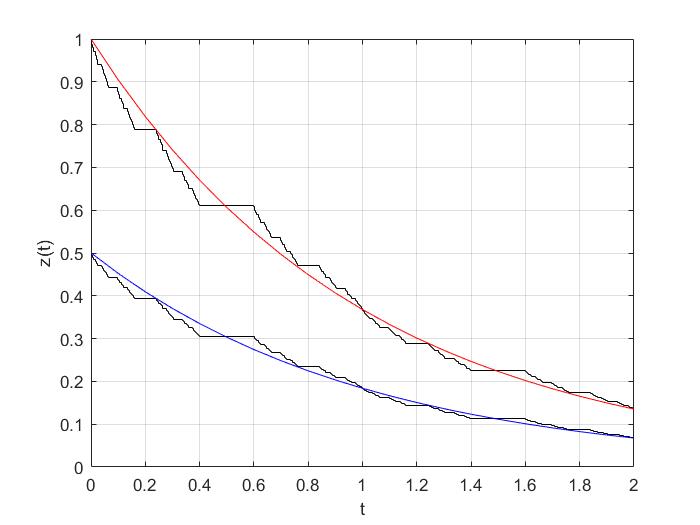}
		\caption{ Solution of Eq. \eqref{luhy}  with $\mu=1/5$ and z(0)=1,0.5 }\label{qwaq1}
\end{subfigure}
	\quad
\begin{subfigure}[t]{0.4\textwidth}
	\centering
	\includegraphics[width=\textwidth]{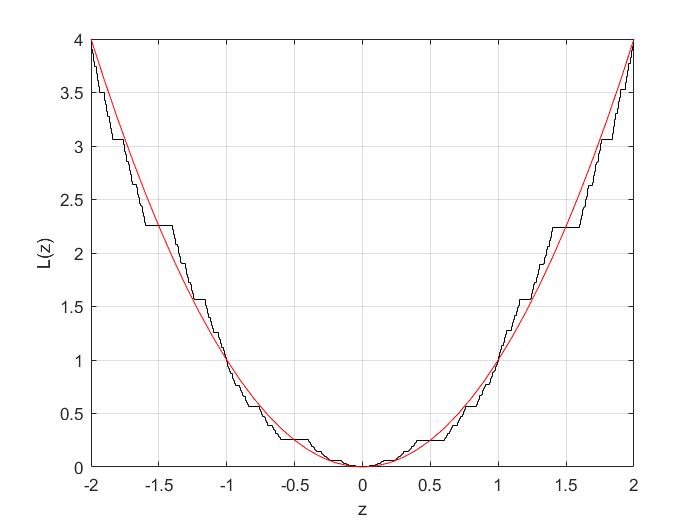}
	\caption{Fractal  Lyapunov function Eq. \eqref{luhy} with $\mu=1/5$}\label{sw2}
	\end{subfigure}
\caption{ Graphs corresponding to  Example 1. }\label{20poidd}
\end{figure}

In this section, we generalize the Lyapunov stability definition for the functions with fractal support.\\
Let us consider the following fractal differential equation with an initial condition
\begin{equation}\label{sqqaq}
  D_{K,t}^{\alpha}h(t)=g(h(t)),~~h(0)=h_{0},
\end{equation}
where $g(h(t)):\Re \rightarrow\Re$, $C^{\mu}=K$ and $g$ has an equilibrium point $h_{e}$, then $g(h_{e})=0$.\\
1) The equilibrium point  $h_{e}$ is called fractal Lyapunov stabile if we have
\begin{equation*}
  \forall~\epsilon>0,~\exists~\delta>0,~~~
  |h(0)-h_{e}|<\delta^{\alpha}\Rightarrow
  |h(t)-h_{e}|<\varepsilon^{\alpha},~~t\geq0.
\end{equation*}
2) The stable equilibrium point $h_{e}$ is said fractal asymptotically stable if
\begin{equation*}
  \forall~\epsilon>0,~\exists~\delta>0,~~~
 |h(0)-h_{e}|<\delta^{\alpha}\Rightarrow
 \lim_{t \rightarrow\infty} |h(t)-h_{e}|= 0.~~
\end{equation*}
3) The equilibrium point $h_{e}$ is called fractal exponentially stable if
\begin{equation*}
  ~\exists~\delta>0,~~~
  |h(0)-h_{e}|<\delta^{\alpha}\Rightarrow
  |h(t)-h_{e}|\leq \kappa^{\alpha}|h(0)-h_{e}|e^{-\lambda \alpha t},~~
\end{equation*}
where $t\geq 0,~\kappa,\lambda\in \Re,$ and $~\kappa>0,~\lambda>0$.\\
Fractal Lyapunov function of  Eq. \eqref{sqqaq} is a function $L:\Re \rightarrow\Re^{+},~R^{+}=[0,+\infty)$ which is  $C^{\mu}$-continuous. Also, its  $\alpha$-order derivative is $C^{\mu}$-continuous. Thus the fractal derivative of $L$ with respect to Eq.\eqref{sqqaq} is written as $L^{*}$ and if it has following condition
\begin{equation}\label{bb12}
 L^{*}=\frac{\partial L}{\partial h}g <0, \forall~ t \in K~\backslash \{0\},
\end{equation}
then, the zero solution of Eq. \eqref{sqqaq} is fractal asymptotically stable.\\
\textbf{Example 1.}~Consider the following fractal differential equation
\begin{equation}\label{luhy}
  D_{K,t}^{\alpha}z(t)=-\chi_{K} z,~~~z(0)=c.
\end{equation}
The general the solution of Eq. \eqref{luhy} is
\begin{equation*}
  z(t)=c\exp(-S_{K}^{\alpha}(t)).
\end{equation*}
A fractal Lyapunov function for studying the stability of Eq. \eqref{luhy} is
 \begin{equation}\label{nxswa}
   L(z)=z^{2}.
 \end{equation}
Then, we have
\begin{equation}\label{uuhy}
  L^{*}=\frac{d L}{dz}(z)=-2z^2<0,~~(z\neq0).
\end{equation}
Thus, the zero solution of Eq.\eqref{luhy} is fractal asymptotically stable (See Figure \ref{qwaq1}).\\
In Figure \ref{20poidd} we have plotted the solution of Eq. \eqref{luhy} in \ref{qwaq1} and corresponding Fractal Lyapunov function in \ref{sw2}.\\
\textbf{Remark.} In Figure \ref{20poidd}, the red curves indicate the orbit  of the solutions of  Eq.\eqref{luhy} and Lyapunov function in the case of  $\alpha=1$.
\section{Qualitative behaviors of solutions of FDEs}
In this section, we present the generalized conditions for the uniform boundedness and convergence of the solutions of the second $\alpha$-order of non-linear  fractal differential equations. On the other hand, we modify and adopt the ordinary calculus conditions in fractal calculus \cite{Cemil-2}. The main results  are obtained using the generalized Lyapunov function with the fractal sets support \cite{Gangal-1,Gangal-2,Alireza-Fernandez-1,Cemil-1,Cemil-2}. \\
Let us consider the following second $\alpha$-order fractal differential equation
\begin{equation}\label{derqas}
  (D_{K,t}^{\alpha})^2y+u(S_{K}^{\alpha}(t))f(y,D_{K,t}
  ^{\alpha}y)D_{K,t}^{\alpha}y+
  v(S_{K}^{\alpha}(t))h(y)
  =q(S_{K}^{\alpha}(t),y,z).~
\end{equation}
By rewritten Eq. \eqref{derqas} in the form of the fractal system  of differential equations and setting $S_{K}^{\alpha}(t)=t'$, we obtain
\begin{eqnarray}\label{121}
  D_{K,t}^{\alpha}y(t') &=& z(t'), \nonumber\\
  D_{K,t}^{\alpha}z(t') &=& -u(t')f(y,D_{K,t}^{\alpha}y)
  z(t')-v(t')h(y)\nonumber\\
  &+&q(t',y,z(t')),
\end{eqnarray}
where $u(t'),~v(t'),~f(y,D_{K,t}^{\alpha}y),~h(y),~z(t')=D_{K,t}^{\alpha}y $ and $q(t',y,z)$ are $C^{\mu}$-continuous functions at every point $t\in C^{\mu}$ and they have well behavior  such that the fractal uniqueness theorem holds for the fractal system \eqref{121}.~Meanwhile, $u(t')$, $v(t')$ are $C^{\mu}$-differentiable on  $C^{\mu}$.\\
\textbf{A. Assumptions}
\begin{enumerate}
\item[(C1)] There are positive constants  $u_{0},~v_{0},~E,~Q$, such that
\begin{eqnarray*}
  1\leq u_{0}^{\alpha}& \leq & u(t')\leq E^{\alpha},\nonumber\\
  1\leq v_{0}^{\alpha}& \leq & v(t')\leq Q^{\alpha},
\end{eqnarray*}
\hspace{-1cm}where we consider $S_{K}^{\alpha}(t)<t^{\alpha}$ \cite{Gangal-1,Gangal-2,Alireza-Fernandez-1}.
  \item[(C2)]   $\lambda_{1}(>0),~\lambda_{2}(>0) \in \Re$ and  $\epsilon_{0},~\epsilon_{1},~\epsilon_{2}$ are small positive numbers such that
\begin{equation*}
        \epsilon_{0}^{\alpha}\leq f(y,D_{K,t}^{\alpha}y).
\end{equation*}
 \item[(C3)] $h(0)=0,~h(y)~\textmd{sgn}(y)>0,~(y\neq0)$, such that
\begin{equation*}
    H(y)=\int_{0}^{y}h(\lambda)d_{K}^{\alpha}
  \lambda\rightarrow\infty~\textmd{as}~|y|\rightarrow\infty,
\end{equation*}
 \hspace{-1cm} and
\begin{equation*}
   0< \lambda_{2}\leq D_{K,t}^{\alpha}h(y).
\end{equation*}
 \item[(C4)]
\begin{equation*}
        \int_{0}^{\infty}\zeta_{0}(t')d_{K}^{\alpha}t<\infty,
        ~D_{K,t}^{\alpha}v(t')\rightarrow0~\textmd{as} ~t\rightarrow\infty,
\end{equation*}
\hspace{-1cm}where $\zeta_{0}(t')=D_{K,t}^{\alpha}v_{+},~D_{K,t}
      ^{\alpha}v_{+}=\max(D_{K,t}^{\alpha}v_{+},0) $.\\

\end{enumerate}
\textbf{Theorem 1.} If assumptions (C1)-(C4)  hold, then the zero solution of Eq.\eqref{derqas} when $q(S_{K}^{\alpha}(t),y,z)=0$ is fractal stable.\\
\textbf{Proof:} For proving this theorem we consider the following fractal Lyapunov function
\begin{equation}\label{Lyggapunov}
  L_{2}(t',y,z)=\int_{0}^{y}h(\lambda)d_{K}^{\alpha}
  \lambda+\frac{z^2}{2v(t')},
\end{equation}
which is positive definite. By calculating fractal time derivative of  \eqref{Lyggapunov}  along the fractal system \eqref{121}, we obtain
\begin{equation*}
  D_{K,t}^{\alpha}L_{2}=-\frac{D_{K,t}^{\alpha}v(t')}{2v(t')^2}z^2
  -\frac{u(t')}{v(t')}f(y,D_{K,t}^{\alpha}y)z^2.
\end{equation*}
We know that $v(t')$ is an increasing function. Hence $ D_{K,t}^{\alpha}v(t')\geq0$. Then, we have
\begin{equation*}
  D_{K,t}^{\alpha}L_{2}\leq -\frac{u(t')}{v(t')}f(y,D_{K,t}^{\alpha}y)z^2 \leq0.
\end{equation*}
In fact, it is obvious that $ L_{2}(t',0,0)=0$ and
\begin{eqnarray*}
  L_{2}(t',y,z)&\geq & \lambda_{2}y^{2}+\frac{z^2}{2v(t')}\nonumber\\
  &\geq & \lambda_{2}y^{2}+\frac{z^2}{2Q^{\alpha}}\nonumber\\
  &\geq &\bar{\lambda}(y^2+z^2),
\end{eqnarray*}
where $\bar{\lambda} =\min (\lambda_{2},\frac{1}{2Q^{\alpha}})$.~Then the proof is complete.\\
\hspace{-1cm}\textbf{B. Assumption}
\begin{enumerate}
\item[(C5)] There are positive constants $0\leq\sigma\leq 1,~0\leq\Delta\leq 1$, such that
\begin{equation*}
\int_{0}^{\infty}r_{i}(t')d_{K}^{\alpha}t,~~~     r_{1}(t'),~r_{2}(t')>0,~~~(i=1,2),
\end{equation*}
\hspace{-1cm}are $C^{\mu}$-Continuous functions and
\begin{equation*}
|q(t',y,z(t'))|\leq r_{1}(t')+r_{2}(t')[H(y)+z^2]^{\sigma'/2}+\Delta^{\alpha}|z|,
\end{equation*}
\hspace{-1cm}where $\sigma'=\sigma^{\alpha}$.
\item[(C6)]
\begin{equation*}
        \epsilon_{0}^{\alpha}\leq f(y,D_{K,t}^{\alpha}y)-\lambda_{1}\leq \epsilon_{1}^{\alpha}.
\end{equation*}
\item[(C7)]
\begin{equation*}
    0\leq\lambda_{2}-D_{K,t}^{\alpha}h(y)\leq \epsilon_{2}^{\alpha}.
\end{equation*}
\end{enumerate}
\textbf{Theorem 2.} Let $q(S_{K}^{\alpha}(t),y,z)\neq0$  and assume $(C1)-(C5)$ hold. Then the solutions of Eq.\eqref{derqas} are fractal uniformly bounded and fractal convergent, namely
\begin{equation}
y(t')\rightarrow0,~D_{K,t}^{\alpha}y(t')\rightarrow 0,~\textmd{as}~ S_{K}^{\alpha}(t)\rightarrow0.
\end{equation}
To prove this theorem, we define a fractal Lyapunov function for Eq.\eqref{derqas} by
\begin{equation}\label{Lyapunov}
  L_{0}(t',y,z)=v(t')\int_{0}^{y}h(\lambda)d_{K}^{\alpha}
  \lambda+\frac{z^2}{2}+k,
\end{equation}
where $k$ is positive constant.\\
Before giving the proof of the above theorem, we present two lemmas, Lemma 1 and Lemma 2, which are needed in the proof of the theorem.\\
\textbf{Lemma 1.} If assumptions  (C1) and (C3) hold, then
\begin{equation*}
E_{1}^{1/\alpha}[H(y)+z^2+k]\leq L_{0}(t',y,z)\leq E_{2}^{1/\alpha}[H(y)+z^2+k],~~ \exists~ E_{1}>0,~E_{2}>0\in\Re.
\end{equation*}
\textbf{Proof:} In view of  assumptions (C1) and (C3) we can derive
\begin{equation*}
L_{0}\geq v_{0}^{\alpha}\int_{0}^{y}h(\lambda)
d_{K}^{\alpha}\lambda+\frac{z^2}{2}+k\geq
E_{1}^{\alpha}[H(y)+z^2+k],
\end{equation*}
where $E_{1}=\min(v_{0},1/2)$.\\
In the same manner, by  assumptions  (C1)  and  (C3),  we can obtain
 \begin{equation*}
L_{0}\leq Q^{\alpha}\int_{0}^{y}h(\lambda)d_{K}
^{\alpha}\lambda+\frac{z^2}{2}+k\leq
E_{2}^{\alpha}[H(y)+z^2+k],
\end{equation*}
where $E_{2}=\max(Q,1)$.~$\Box$\\
\textbf{Lemma 2.} If assumptions   (C1)-(C4) are valid, then
\begin{equation*}
 \exists~E_{3}>0,~E_{4}>0,~D_{K,t}^{\alpha}L_{0}\leq-
 E_{3}^{\alpha}z^2+(r_{1}(t')+r_{2}(t'))|z|+
r_{2}(t')[H(y)+z^2]+E_{4}^{\alpha}\zeta_{0}L_{0},
\end{equation*}
where $E_{3},E_{4}\in\Re$.\\
\textbf{Proof:} Fractal differentiating of the fractal Lyapunov function \eqref{Lyapunov} along with fractal system  \eqref{121},  we get
\begin{equation*}
  D_{K,t}^{\alpha}L_{0}=-u(t')f(y,z)z^2+q(t',y,z)z+D_{K,t}^{\alpha}v(t')
  \int_{0}^{y}h(\lambda)d_{K}^{\alpha}\lambda.
\end{equation*}
By using the assumptions of the Theorem 2, we obtain
\begin{eqnarray*}
  D_{K,t}^{\alpha}L_{0} &\leq& -E^{\alpha}(\lambda_{1}+\epsilon_{0})z^2+|q(t',y,z)||z|+
  \zeta_{0}H(y)\nonumber \\
   &\leq& -2E_{3}^{\alpha}z^2+|q(t',y,z)||z|+E_{4}^{\alpha}\zeta_{0}L_{0},
\end{eqnarray*}
where
\begin{equation*}
  E_{3}=E(\lambda_{1}+\epsilon_{0})/2,~E_{4}=(1/E_{1}).
\end{equation*}
Here, in view of (C5), the term $|q(t',y,z)||z|$ can be written as
\begin{equation*}
  |q(t',y,z)||z|\leq\bigg(r_{1}(t')+r_{2}(t')[H(x)+z^2]^
  {\sigma'/2}\bigg)|z|+\Delta^{\alpha} z^2.
\end{equation*}
Hence, we have
\begin{equation*}
  D_{K,t}^{\alpha}L_{0} \leq-2E_{3}^{\alpha}z^2+(r_{1}(t')+r_{2}(t')[H(x)+z^2]^{\sigma'/2})
  |z|+\Delta^{\alpha} z^2+E_{4}^{\alpha}\zeta_{0}L_{0}.
\end{equation*}
Set $\Delta=E_{3}$. Then
\begin{equation}\label{gtb}
  D_{K,t}^{\alpha}L_{0} \leq-E_{3}^{\alpha}z^2+(r_{1}(t')+r_{2}(t')[H(x)+z^2]^{\sigma'/2})
  |z|+E_{4}^{\alpha}\zeta_{0}L_{0}.
\end{equation}
Using the following inequality
\begin{equation}\label{kkgtb2}
  [H(x)+z^2]^{\sigma'/2}\leq 1+[H(x)+z^2]^{1/2},
\end{equation}
 taking into account \eqref{gtb} and  \eqref{kkgtb2} we obtain
\begin{equation*}
  D_{K,t}^{\alpha}L_{0} \leq-E_{3}^{\alpha}z^2+(r_{1}(t')+r_{2}(t'))|z|+
  r_{2}(t')[H(x)+z^2]+
  E_{4}^{\alpha}\zeta_{0}L_{0}.~
\end{equation*}
To complete the proof of the theorem, we consider the fractal Lyapunov function  $L_{0}$ defined by
\begin{equation}\label{55847}
  L(t',y,z)=e^{-\int_{0}^{t}\zeta(\theta)
  d_{K}^{\alpha}\theta}L_{0}(t',y,z),
\end{equation}
where
\begin{equation*}
  \zeta(t')=E_{4}^{\alpha}\zeta_{0}+\frac{4}
  {E_{1}^{\alpha}}(r_{1}(t')+r_{2}(t')),
\end{equation*}
and
\begin{equation}\label{mmmm}
  \psi_{1}(||\bar{y}||)\leq V(t',y,z)\leq \psi_{2}(||\bar{y}||),
\end{equation}
with $\bar{y}=(y,z)\in \Re^2$ and $t\in C^{\mu}$, and $\psi_{1},\psi_{2}$ are $C^{\mu}$-continuous and increasing functions such that $\psi_{1}(||\bar{y}||)\rightarrow\infty$ while $||\bar{y}||\rightarrow\infty$.\\
Fractal differentiating fractal Lyapunov function \eqref{55847} and considering  the fractal system \eqref{121} and assumptions of  Theorem 2, we have
\begin{align*}
 D_{K,t}^{\alpha} L(t',y,z) & =e^{-\int_{0}^{t}\zeta(\theta)
  d_{C^{\mu}}^{\alpha}\theta}[ D_{K,t}^{\alpha}L_{0}-\zeta(t')L_{0}]\nonumber\\
   & \leq e^{-\int_{0}^{t}\zeta(\theta)d_{K}^{\alpha}\theta}[-E_{3}^
   {\alpha}
   z^2+(r_{1}(t')+r_{2}(t'))|z|+
   r_{2}(t')[H(y)+z^2]\nonumber\\&-4(r_{1}(t')+r_{2}(t'))(H(y)+z^2+k)]
  \nonumber\\
   & \leq e^{-\int_{0}^{t}\zeta(\theta)
  d_{K}^{\alpha}\theta}[-E_{3}^{\alpha}z^2+(r_{1}(t')+r_{2}(t'))|z|
  \nonumber\\&-2 (r_{1}(t')+r_{2}(t'))(H(y)+z^2+k)] \nonumber\\
   & \leq e^{-\int_{0}^{t}\zeta(\theta)
  d_{K}^{\alpha}\theta}[-E_{3}^{\alpha}z^2-2(r_{1}
  (t')+r_{2}(t'))\{(|z|-
  \frac{1}{4})^2-\frac{1}{16}+2k\}].
\end{align*}
If we choose $k\geq\frac{1}{32}$, then it follows that there exists a positive $E_{5}$ such that
\begin{equation}\label{l59n}
  D_{K,t}^{\alpha} L(t',y,z) \leq -E_{5}^{\alpha}z^2.
\end{equation}
By considering \eqref{mmmm} and \eqref{l59n} it follows that  all solutions of Eq.\eqref{121} are fractal uniformly bounded.\\
Consider the fractal system of  differential equations
\begin{equation}\label{ppp}
 D_{K,t}^{\alpha} \bar{y}=M(t',\bar{y})+N(t',\bar{y}),
\end{equation}
where $M,N$ are $C^{\mu}\subset \Re^{+}$-continuous and vector functions, and $C^{\mu}\times Q\subset \Re^2$  is an open set. Moreover, it is clear that
\begin{equation*}
  ||N(t',\bar{y})||\leq N_{1}(t',\bar{y})+N_{2}(\bar{y}),
\end{equation*}
where $N_{1}(t',\bar{y}),N_{2}(\bar{y})$ are $C^{\mu}$-continuous and non-negative functions. $\Box$\\
\textbf{Lemma 3.} Let $L: C^{\mu}\times Q$ be a function $C^{\mu}$-continuous and $C^{\mu}$-differentiable such that
\begin{equation*}
  D_{K,t}^{\alpha}L(t',\bar{y})\leq-B(||\bar{y}||),
\end{equation*}
where $B(||\bar{y}||)$ is a positive  definite in the closed set $\Psi\in Q$ and $M(t',\bar{y})$ satisfies the following.
\begin{enumerate}
  \item $M(t',\bar{y})$ tends to $K(\bar{y})$ for $\bar{y}\in \Psi$ as $t\rightarrow \infty$, $K(\bar{y})$ is a $C^{\alpha}$-continuous on $\Psi$.
  \item
  \begin{equation*}
    \forall~ \epsilon>0,~ \bar{y}\in \Psi,~\exists~ \delta=\delta(\epsilon, \bar{z})>0,~ T=T(\epsilon,~ \bar{z})>0,
  \end{equation*}
  such that if
  \begin{equation*}
    t\geq T,~||\bar{y}-\bar{z}||<\delta(\epsilon,\bar{z}),
  \end{equation*}
  we have
  \begin{equation*}
    ||M(t',\bar{y})-M(t'-\bar{z})||<\epsilon^{\alpha}.
  \end{equation*}
  \item $N_{2}(\bar{y})$ is positive definite on closed  $\Psi$ of $Q$.
\end{enumerate}
Then every bounded solution of Eq.\eqref{ppp} approaches to the fractal system
\begin{equation}
 D_{K,t}^{\alpha} \bar{y}=\mathfrak{K}(\bar{y}),
\end{equation}
which is contained in $\Psi$ as $t\rightarrow \infty.$\\
\textbf{Proof.} Now, we consider \eqref{121}. It follows that
\begin{equation*}
  M(t',\bar{y})=\left(
                  \begin{array}{c}
                    z \\
                    -u(t')f(y,z)z-v(t')h(y) \\
                  \end{array}
                \right)
\end{equation*}
and
\begin{equation*}
  N(t',\bar{y})=\left(
                  \begin{array}{c}
                    0 \\
                    q(t',y,z) \\
                  \end{array}
                \right).
\end{equation*}
Then
\begin{equation*}
  ||N(t',\bar{y})||\leq(r_{1}(t')+r_{2}(t'))
  [H(y)+z^2]^{\sigma'/2}+\Delta^{\alpha}|z|.
\end{equation*}
We can also write
\begin{equation*}
  N_{1}(t',\bar{y})=r_{1}(t')+r_{2}(t')[H(y)+z^2]^{\sigma'/2}
\end{equation*}
and
\begin{equation*}
  N_{2}(\bar{y})=\Delta^{\alpha}|z|.
\end{equation*}
The functions $M(t',\bar{y})$ and $N(t',\bar{y})$ satisfy the conditions of Lemma 3. Set $\psi_{1}(||\bar{y}||)=E_{5}^{\alpha}z^2$, then
\begin{equation*}
  D_{F,t}^{\alpha}L(t',y,z)\leq-\psi_{1}(||\bar{y}||),
\end{equation*}
where the function $||\bar{y}||$ is positive definite on $\Psi=\{(y,z)|y\in \Re,z=0\}$. We get
\begin{equation*}
  M(t',\bar{y})=\left(
                  \begin{array}{c}
                    0 \\
                    -v(t')h(y) \\
                  \end{array}
                \right)
\end{equation*}
by using  (C4) condition  of Theorem 1.~If  we suppose
\begin{equation}\label{50}
  \mathfrak{K}(\bar{y})=\left(
               \begin{array}{c}
                 0 \\
                 -v_{\infty}h(y) \\
               \end{array}
             \right),
\end{equation}
then all the conditions of Lemma 3 are satisfied. It is straight forward to see that $N_{2}(\bar{y})$ is  positive  definite function. Since the solutions of fractal system \eqref{121} are bounded, therefore by using Lemma 3 we have
\begin{equation*}
  D_{K,t}^{\alpha} \bar{y}=\mathfrak{K}(\bar{y}),
\end{equation*}
which is semi-invariant set of the fractal system contained in $\Psi$ as $t\rightarrow \infty$. In view of \eqref{50}, we have following
\begin{equation}\label{jjj}
  D_{K,t}^{\alpha} y=0,~~~D_{K,t}^{\alpha}z=-v_{\infty}h(y).
\end{equation}
Fractal system Eq.\eqref{jjj} has solution
\begin{equation*}
  y=c_{1},~~~z=c_{2}-v_{\infty}h(c_{1})
  (S_{K}^{\alpha}(t)-S_{K}^{\alpha}(t_{0})).
\end{equation*}
In order to remain in $\Psi$, the solutions must be
\begin{equation*}
  c_{2}-v_{\infty}h(c_{1})
  (S_{K}^{\alpha}(t)-S_{K}^{\alpha}(t_{0}))=0, \forall ~t\geq t_{0},
\end{equation*}
which implies $k=0,h(c_{1})=0$, so that $c_{1}=c_{2}=0$. Then $\bar{y}=\bar{0}$ is the solution of $D_{K,t}^{\alpha} \bar{y}=\mathfrak{K}(\bar{y})$ remaining in $\Psi$. Consequently, we arrive at
\begin{equation*}
  y(t')\rightarrow0,~~~D_{K,t}^{\alpha} y(t')\rightarrow0,~as~t\rightarrow0.~\Box
\end{equation*}
\textbf{Example 2.} Consider the fractal differential equation
\begin{equation}\label{yytrewq}
  (D_{K,t}^{\alpha})^2y(t)+s(y)D_{K,t}^{\alpha}y(t)+h(y)=0.
\end{equation}
This is equivalent to the fractal system
\begin{eqnarray}\label{1ik}
  D_{K,t}^{\alpha}y&=&z-S(y)\nonumber\\
  D_{K,t}^{\alpha}z&=&-h(y),
\end{eqnarray}
where
\begin{equation*}
  S(y)=\int_{0}^{y}s(\beta)d_{K}^{\alpha}\beta.
\end{equation*}
and, if we suppose
\begin{equation*}
  H(y)=\int_{0}^{y}h(\rho)d_{K}^{\alpha}\rho,
\end{equation*}
Consider the fractal function
\begin{equation}\label{147}
  L_{1}(y,z)=H(y)+\frac{z^2}{2},
\end{equation}
which is a strong fractal Lyapunov function. Fractal differentiating of  $L_{1}(y,z)$ and Eq.\eqref{147} we get
\begin{equation*}
  D_{K,t}^{\alpha}L_{1}(y,z)=
  h(y)D_{K,t}^{\alpha}y+zD_{K,t}^{\alpha}z=-h(y)S(y).
\end{equation*}
Using the assumption theorem we obtain
\begin{equation*}
  D_{K,t}^{\alpha}L_{1}(y,z)=-h(y)S(y)<0,
\end{equation*}
which shows that the solutions of Eq.\eqref{yytrewq} are fractal uniformly bounded and fractal ultimately bounded.\\
\textbf{Example 3.} Consider harmonic oscillator on the fractal time as follows:
\begin{equation}\label{zswesaq}
  (D_{K,t}^{\alpha})^2y(t)+\mathfrak{C}_{K}y(t)=0,~t\in K,~\mathfrak{C}_{K}>0,
\end{equation}
where $\mathfrak{C}_{K}$ is constant. The equivalent fractal system is
\begin{eqnarray}
  D_{K,t}^{\alpha}y &=& z, \nonumber\\
  D_{K,t}^{\alpha}z &=& - \mathfrak{C}_{K}y.
\end{eqnarray}
The fractal Lyapunov function correspond to Eq.\eqref{zswesaq} is
\begin{equation}\label{nm86}
  L(y,z)=\frac{1}{2}\mathfrak{C}_{K}y^2+\frac{1}{2}z^2,
\end{equation}
where $L(0,0)=0$ and $L(y,z)>0$ for $(y,z)\in \Re^{2}\backslash(0,0)$.
\begin{figure}[H]
\centering	
\begin{subfigure}[t]{0.4\textwidth}
		\centering
		\includegraphics[width=\textwidth]{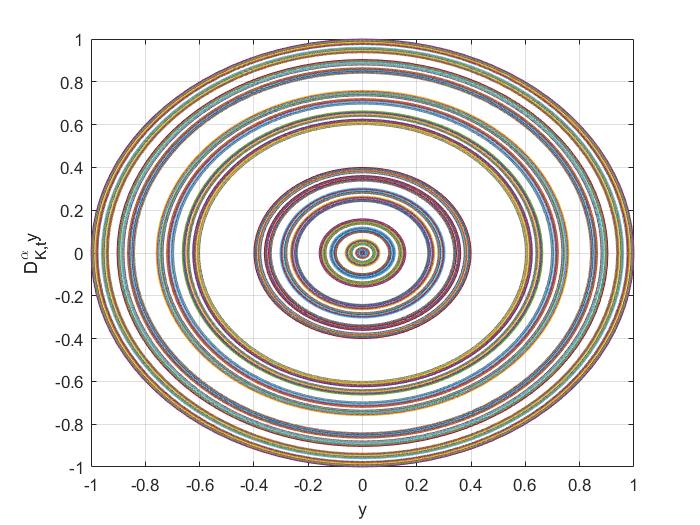}
		\caption{ Solution of Eq. \eqref{zswesaq}  with $\mu=1/5$ and }\label{qwazaq1}
\end{subfigure}
	\quad
\begin{subfigure}[t]{0.4\textwidth}
	\centering
	\includegraphics[width=\textwidth]{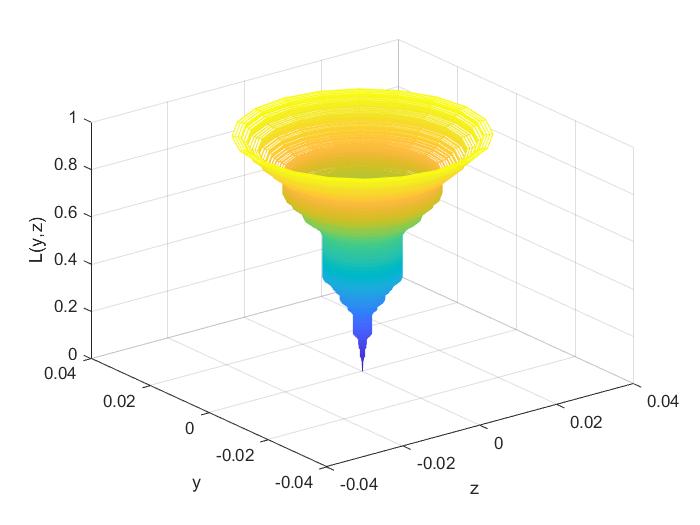}
	\caption{Fractal  Lyapunov function Eq. \eqref{zswesaq} with $\mu=1/5$}\label{22}
	\end{subfigure}
\caption{Graphs corresponding to  Example 2. }\label{99poi}
\end{figure}
 Then, it is obtain that
\begin{equation*}
  D_{K,t}^{\alpha}L(y,z)=\frac{\partial}{\partial y}L(y,z)D_{K,t}^{\alpha}y+\frac{\partial}{\partial y},
  L(y,z)D_{K,t}^{\alpha}z=0.
\end{equation*}
Hence, the zero solution $(0,0)$ is a fractal stable point. In Figure \ref{99poi}, we have sketched solutions of  Eq. \eqref{zswesaq} and Eq. \eqref{nm86}.

\section{Conclusion}
In this paper, we have suggested conditions for the fractal stability, uniformly boundedness and the asymptotic behaviors of solutions of second $\alpha$-order fractal differential equations. The analogous theorems of stability, uniformly boundedness and asymptotic behavior to standard calculus have been given and adopted in fractal calculus. The generalized conditions include solutions and functions which are non-differentiable in sense of ordinary calculus.

\textbf{Acknowledgement:} \\
This research was completed with the support of the Scientific and Technological Research Council of Turkey (T\"{U}B\.{I}TAK) (2221-Fellowships for Visiting Scientists and Scientists on Sabbatical Leave – 2221-2018/3 period) when Alireza Khalili Golmankhaneh was a visiting scholar at Van Yuzun\c{c}u Yil University, Van, Turkey

\end{document}